\documentclass{amsart}

\usepackage{amssymb}

\def\bc{\begin{center}}
\def\ec{\end{center}}
\def\be{\begin{equation}}
\def\ee{\end{equation}}
\def\ben{\begin{enumerate}}
\def\een{\end{enumerate}}
\def\bfg{\begin{figure}}
\def\efg{\end{figure}}
\def\bq{\begin{quote}}
\def\eq{\end{quote}}
\def\bd{\begin{description}}
\def\ed{\end{description}}

\def\h{\hbar}
\def\p{\partial}
\def\w{\wedge}

\def\dim{\operatorname{dim}}

\newcommand{\CC}{{\mathbb C}}

\newcommand{\ZZ}{{\mathbb Z}}
\newcommand{\QQ}{{\mathbb Q}}

\newcommand{\lan}{\langle}
\newcommand{\ran}{\rangle}
\newcommand{\Gb}{{\mathfrak B}}

\newcommand{\gl}{\lambda}

\newcommand{\go}{\omega}

\newcommand{\gs}{\sigma}

\newcommand{\gO}{\Omega}
\newcommand{\gS}{\Sigma}

\newcommand{\J}{{\mathcal J}}
\newcommand{\K}{{\mathcal K}}

\newcommand{\Tau}{{\mathcal T}}

\newcommand{\M}{\overline{\mathcal M}}
\newcommand{\ct}{\operatorname{ct}}
\newcommand{\ev}{\operatorname{ev}}
\newcommand{\ft}{\operatorname{ft}}
\renewcommand{\QQ}{\mathbb Q}
\renewcommand{\CC}{\mathbb C}
\renewcommand{\t}{\mathbf t}
\newcommand{\q}{\mathbf q}
\newcommand{\f}{\mathbf f}
\newcommand{\g}{\mathbf g}

\newcommand{\x}{\mathbf x}

\newcommand{\D}{\mathcal D}

\newcommand{\A}{\mathcal A}

\renewcommand{\a}{\alpha}
\renewcommand{\b}{\beta}
\renewcommand{\c}{\gamma}
\renewcommand{\d}{\delta}
\renewcommand{\H}{{\mathcal H}}
\newcommand{\F}{{\mathcal F}}
\newcommand{\G}{{\mathcal G}}

\renewcommand{\L}{{\mathcal L}}

\newcommand{\1}{{\bf 1}}
\newcommand{\s}{{\mathbf s}}

\renewcommand{\l}{{\lambda}}
\newcommand{\ch}{\operatorname{ch}}

\newcommand{\e}{{\mathbf e}}
\newcommand{\td}{\operatorname{td}}
\newcommand{\0}{{\mathbf 0}}
\newcommand{\Res}{\operatorname{Res}}

\newcommand{\Td}{\operatorname{Td}}
\newcommand{\Ch}{\operatorname{Ch}}
\newcommand{\qCh}{\operatorname{qCh}}
\newcommand{\qch}{\operatorname{qch}}

\begin{document}

\title[Symplectic geometry of Frobenius structures]
{Symplectic geometry of Frobenius structures}
\author{Alexander B. Givental} 
\address{UC Berkeley} 
\thanks{Research is partially supported by NSF Grant DMS-0072658} 


\maketitle

The concept of a Frobenius manifold was introduced by B. Dubrovin \cite{D}
to capture in an axiomatic form the properties of correlators 
found by physicists (see \cite{DW}) in two-dimensional topological 
field theories ``coupled to gravity at the tree level''.
The purpose of these notes is
to reiterate and expand the viewpoint, outlined in the paper \cite{CG} 
of T. Coates and the author, which recasts this concept
in terms of linear symplectic geometry and exposes the    
role of the twisted loop group $\L^{(2)}GL_N$ of hidden symmetries.

We try to keep the text introductory and non-technical. 
In particular, we supply details of some simple results 
from the axiomatic theory, including a several-line proof of
the genus $0$ Virasoro constraints not mentioned elsewhere, 
but merely quote and refer to the literature for a number of 
less trivial applications, such as the quantum 
Hirzebruch--Riemann--Roch theorem in the theory 
of cobordism-valued Gromov--Witten invariants.
The latter is our joint work in progress with Tom Coates,
and we would like to thank him for numerous discussions of the
subject. 

The author is also thankful to Claus Hertling, Yuri Manin and
Matilde Marcolli for the invitation to the workshop on 
Frobenius structures held at MPIM Bonn in Summer 2002. 

\medskip

{\bf Gravitational descendents.}  
In Witten's formulation of topological gravity \cite{W} 
one is concerned with certain ``correlators''
called {\em gravitational descendents}. The totality of genus $0$ 
gravitational descendents is organized as the 
set of Taylor coefficients of a single formal function $\F$ of a sequence 
of vector variables $\t = (t_0, t_1, t_2, ... )$. 
The vectors $t_i$ are elements of a finite-dimensional vector 
space\footnote{or super-space, but we will systematically ignore the signs 
that may come out of this possibility.}  
which we will denote $H$. One assumes that $H$ is equipped 
with a symmetric 
non-degenerate bilinear form $(\cdot,\cdot)$ and with a distinguished  
non-zero element $\1 $. 
The axioms are formulated as certain partial differential 
equations on $\F$: an infinite set of {\em Topological Recursion Relations} 
(TRR), the {\em String Equation} (SE) and the {\em Dilaton Equation} (DE). 
To state the axioms using the following coordinate notation: 
introduce a basis $\{ \phi_{\a} \} $ in $H$ such that 
$\phi_1=\1$, put $t_k=\sum t_k^{\a} \phi_{\a}$, 
$(\phi_{\a},\phi_{\b}) = g_{\a\b}$, $(g^{\a\b}) = (g_{\a\b})^{-1}$ and write 
$\p_{\a,k}$ for the partial derivative $\p/\p t_k^{\a}$.
Then the axioms read:
\begin{equation} \tag{DE} \p_{1,1} \F\ (\t)  = 
\sum_{n=0}^{\infty} \sum_{\nu} t^{\nu}_n \p_{\nu,n} \F\ (\t)\ -\ 2\F\ (\t),
\end{equation}  
\begin{equation} \label{(SE)} \tag{SE}  \p_{1,0} \F \ (\t) = 
\frac{1}{2}(t_0,t_0)+\sum_{n=0}^{\infty} 
\sum_{\nu} t_{n+1}^{\nu} \p_{\nu,n} \F\ (\t), \end{equation} 
\begin{equation} \label{TRR} \tag{TRR}   
 \p_{\a,k+1}\p_{\b,l}\p_{\c,m} \F\ (\t) = \sum_{\mu\nu}
 \p_{\a,k} \p_{\mu,0} \F\ (\t) \ g^{\mu\nu} \ 
\p_{\nu,0}\p_{\b,l}\p_{\c, m}\F\ (\t) \end{equation}  
for all $\a,\b,\c$ and any $k,l,m\geq 0$.  
We will see soon that the system DE$+$SE$+$TRR admits a transparent 
geometrical interpretation.

Let $\H$ denote the {\em loop space} $H((z^{-1}))$ consisting of Laurent 
series in $1/z$ with vector coefficients. 
Introduce the symplectic form $\gO$ on $\H$:
\[ \gO (\f,\g )  = \frac{1}{2\pi i} \oint\ ( \f (-z), \g(z) ) \ dz .\]
The polarization $\H = \H_{+}\oplus \H_{-}$ by the Lagrangian subspaces
$\H_{+} = H[z]$ and $\H_{-} = z^{-1} H [[ z^{-1} ]]$ provides a symplectic
identification of $(\H, \gO)$ with the cotangent bundle $T^*\H_{+}$. 
We encode the domain variables $\t = (t_0,t_1,t_2,...)$ of the 
function $\F$ in the coefficients of the vector polynomial
\[ \q (z) = - z+ t_0+t_1z+t_2z^2 + ...  \qquad 
\text{(where $-z := -\1 z$)}\ .\]  
Thus the formal function $\F (\t)$ near $\t=0$ becomes a formal function 
on the space $\H_{+}$ near the shifted origin $\q = -z$. This convention is 
called the {\em dilaton shift}. Denote by $\L$ the graph of 
the differential $d\F$:
\[ \L = \{ ({\mathbf p}, \q) \in T^*\H_{+}:\ {\mathbf p} = d_{\q} \F\} .\]
This is a formal germ at $\q = -z$ of a Lagrangian section of 
the cotangent bundle $T^*\H_{+}$ and can therefore be considered as
a formal germ of a Lagrangian submanifold in the symplectic loop space 
$(\H,\gO)$.

\medskip

{\bf Theorem 1.} {\em The function $\F$ satisfies 
TRR, SE and DE if and only if the corresponding 
Lagrangian submanifold $\L \subset \H$ has the following property:
\begin{equation} \tag{$\star$} \begin{array}{l} 
\text{$\L$ is a Lagrangian cone with the vertex at the origin and such} \\ 
\text{that its tangent spaces $L$ are tangent to $\L$ exactly along $zL$.}
\end{array} \end{equation}}

\medskip

More precisely, the last condition is interpreted as follows. First, it
shows that the tangent spaces $L=T_{\f}\L$ satisfy $zL \subset L$
(and therefore $\dim L/zL = \dim \H_{+} /z\H_{+} =\dim H$). Second,
it implies that $zL \subset \L$. Third, the same $L$ is the tangent 
space to $\L$ not only along the line of $\f$ but also at all smooth 
points in $zL$. This is meant to include $\f$, i.e. $L=T_{\f}\L$ implies
$\f \in zL$.


The tangent spaces $L$ vary therefore in a $\dim H$-parametric family. 
In particular they generate a {\em variation of    
semi-infinite Hodge structures} in the sense of S. Barannikov \cite{B},
i.e. a family of semi-infinite flags $ \cdots zL \subset L \subset z^{-1}L 
\cdots$ satisfying the Griffiths integrability condition.

Thus another way to phrase the property ($\star$) is to say that
 
\noindent 
{\em $\L$ is a cone ruled by the isotropic subspaces $zL$ varying in a 
$\dim H$-parametric family with the tangent spaces along $zL$ equal to
the same Lagrangian space $L$.} 

Finally, the property should be
understood in the sense of formal geometry near a point $\f_0$ 
different from the origin. In particular, the fact that cones develop 
singularities at the origin does not contradict the fact that $\L$ has 
Lagrangian tangent spaces as a section of $T^*\H_{+}$ near $\q=-z$. 

\medskip

{\em Proof.} Let us begin with the cone $\L$ satisfying ($\star$).
Using the polarization $\H = T^*\H_{+}$ we recover $\F$ 
as the restriction of the quadratic form ${\mathbf p} (\q)/2$ 
(evaluation of the covector on a vector) to the cone $\L$ parameterized by
$\H_{+}$ as a section of the cotangent bundle. The resulting function
$\F$ is homogeneous of degree $2$ (since $\L$ is a cone) which is 
exactly what DE expresses after the dilaton shift. 

Multiplication by $z^{-1}$ in $\H$ is anti-symmetric with respect 
to the symplectic form $\gO$ and thus defines a linear hamiltonian vector
field $\f \mapsto z^{-1}\f$ corresponding to the quadratic hamiltonian
$\gO (z^{-1} \f, \f)/2$. The vector field is tangent to $\L$ since
$z^{-1}\f$ is contained in the tangent space $L$ at the point $\f \in zL$,
The cone $\L$ is therefore situated in the zero level of the quadratic 
hamiltonian. It is straightforward to check that this condition
written as the Hamilton--Jacobi equation for the generating function $\F$
coincides with the string equation.

Next, the property of $\L$ to have constant tangent spaces $L$ 
along subspaces $zL$ of codimension $\dim H$ shows that the quadratic
differentials $d^2_{\t} \F$ depend only on $\dim H$ variables.
In particular this applies to $\tau^{\d}:=\sum g^{\d\mu}\p_{\mu,0}\p_{1,0}\F$.
Differentiating the string equation in $t^{\mu}_0$ we find 
\begin{equation}  \tau^{\d} = t^{\d}_0+\sum_{n,\mu,\nu} 
g^{\d\mu}\ t^{\nu}_{n+1}\ \p_{\nu,n}\p_{\mu,0}\F. \end{equation}
In particular
$\tau^{\d} (\t)|_{t_1=t_2=...=0} = t^{\d}_0$ are independent
and thus $\tau^{\d}$ can be taken in the role of the $\dim H$ variables.

Note that the tangent spaces $L$ are transverse to $\H_{-}$ and so the spaces
$zL$ are transverse to $z\H_{-}$. Therefore the projection of $zL$ to the base
$\H_{+}$ along $\H_{-}$ is transverse to the slice $q_1=-\1, q_2=q_3=...=0$,
and the intersection with the slice has $q_0=\tau$. Making use of 
the correlator notation
\[  \lan \phi_{\a}\psi^k, \psi_{\b}\psi^l, ..., \phi_{\c}\psi^m \ran (\tau) := 
 (\p_{\a,k}\p_{\b,l}...\p_{\c,m} \F) | \ t_0=\tau, t_1=t_2=...=0, \]
we can write
\begin{equation}  
\p_{\b,l}\p_{\c,m} \F (\t ) = \lan \phi_{\b}\psi^l, \phi_{\c}\psi^m \ran 
(\tau (\t)) . \end{equation} 

Consider (1) as an implicit equation for $\tau(\t)$:
\begin{equation} G^{\d}(\tau, \t):= \tau^{\d} - t^{\d}_0-\sum_{n,\mu,\nu} 
g^{\d\mu}\ t^{\nu}_{n+1}\ \lan \phi_{\nu}\psi^l, \phi_{\mu} \ran (\tau) = 0. 
\end{equation}
Using the rules of implicit differentiation we conclude that the
matrix $(\p_{\l,0} \tau^{\d})$ is inverse to $(\p G^{\d}/\p \tau^{\l})$ at
$\tau=\tau(\t)$. Using this together with implicit differentiation once again,
we find
\[ \p_{\a,k+1} \tau^{\d} = \sum_{\mu,\nu} \p_{\nu,0}\tau^{\d}\ g^{\mu\nu}\ 
\lan \phi_{\a}\psi^k, \phi_{\mu}\ran (\tau(\t)).\]
This is a special case of TRR whose general form is obtained now by
differentiation of (2):
\begin{align} \p_{\a,k+1}\p_{\b,l}\p_{\c,m}\F & =  
 \sum_{\d}  
\p_{\a,k+1}\tau^{\d}\ \lan \phi_{\d}, \phi_{\b}\psi^l,\phi_{\c}\psi^m\ran 
\\ & =  \sum_{\d,\mu,\nu } \lan \phi_{\a}\psi^k, \phi_{\mu}\ran \ g^{\mu\nu}\ 
\p_{\nu,0} \tau^{\d} \ 
\lan \phi_{\d}, \phi_{\b}\psi^l,\phi_{\c}\psi^m\ran 
\\ & = \sum_{\mu,\nu} \lan \phi_{\a}\psi^k, \phi_{\mu}\ran\ g^{\mu\nu}\ 
\p_{\nu,0}\p_{\b,l}\p_{\c,m} \F .\end{align}  

\medskip    

In the opposite direction, suppose that $\F$ satisfies TRR, SE and DE. 
It is straightforward to see that given 
$\lan \phi_{\b}\psi^l,\phi_{\c}\psi^l\ran (\tau)$,
the TRR uniquely determines all the other correlators with at least two seats.

We use (3) to implicitly define the (formal) map $\t \mapsto \tau(\t)$ and
try the right hand side of (2) in the role of $\p_{\b,l}\p_{\c,m}\F$.
Repeating the above derivation (4,5,6) we conclude that TRR is satisfied
and hence, by the uniqueness, (2) holds true indeed.  

We see therefore that the tangent spaces
to the graph $\L$ of $d\F$ are constant along the fibers of the map
$\t \mapsto \tau$. 
[ These fibers are affine subspaces (of codimension $\dim H$)
since $G$ is linear inhomogeneous in $\t$ and $\tau (t_0,0,0,...) = t_0$.
Taking $\t= (0,1,0,0,...)$ in $G$ we find 
$G^{\d}= \tau^{\d} -\sum g^{\d,\mu} \lan \phi_{\mu}, \1 \ran (\tau)$. 
This equals $0$ for all $\tau$ because $\lan \phi_{\mu}, \1 \ran (\tau) =
(\phi_{\mu} ,\tau )$ as it follows from the string equation. Thus the fibers
of the map $\t \mapsto \tau $ become linear subspaces after the dilaton shift.]
 
Due to the dilaton equation, $\F$ as a function on $\H_{+}$ is homogeneous of
degree $2$, and hence $\L$ is homogeneous of degree $1$. Thus $\L$ is a cone,
so the tangent spaces $L = T_{\f}\L$ contain the application points $\f$.
Recalling that the string equation expresses the tangency of $z^{-1}\f$ to
$\L$ at the points $\f \in \L$ we see that in fact $\f \in L\cap zL$. 
This is true for all those $\f$ where the tangent space is $L$. Therefore 
the projection of $zL$ to $\H_{+}$ along $\H_{-}$ contains one
of the fibers of the map $\t \mapsto \tau$. In fact the projection
coincides with the fiber since both have codimension $\dim H$ in $\H_{+}$. 
For the space $zL$ itself this means that it coincides with the set of points 
$\f\in \L$ where $T_{\f}\L=L$. In particular $zL\subset L$. $\square $

\medskip

The property ($\star$) 
is formulated in terms of the symplectic structure $\gO$
and the operator of multiplication by $z$ but it does not depend on the 
choice of the polarization. This shows that the system DE$+$SE$+$TRR
has a huge group of hidden symmetries. Let $\L^{(2)}GL(H)$ denote the 
{\em twisted loop group} which consists of $\operatorname{End}(H)$-valued 
formal Laurent series $M$ in the indeterminate $z^{-1}$ satisfying 
$M^*(-z)M(z)=\1$. Here $\ ^*$ denotes the adjoint with respect to 
$(\cdot ,\cdot )$. Later we will mention a number of advantages that 
the following corollary has to offer. 

\medskip

{\bf Corollary.} {\em The action of the twisted loop group 
preserves the class of the Lagrangian cones $\L$
satisfying ($\star$) and, generally speaking, yields 
new generating functions $\F$ which satisfy the system DE$+$SE$+$TRR whenever
the original one does.}
        
\medskip

{\bf Frobenius structures.}
By a result of S. Barannikov \cite{B}, 
the set of tangent spaces $L$ to the 
cone $\L$ satisfying ($\star$) carries a canonical Frobenius 
structure. We quot his argument and outline the construction. 

Consider the intersection of the cone $\L$ with the affine space
$-z + z\H_{-}$. The intersection is parameterized by $\tau\in H$ via
its projection to $-z+H$ along $\H_{-}$ and can be considered as the graph
of a function from $H$ to $\H$ called the {\em $J$-function}:
\[  \tau  \mapsto J(-z,\tau) = -z + \tau + \sum_{k>0} J_k (\tau) (-z)^{-k}. \]
The derivatives $\p J/\p \tau^{\d} = \phi_{\d} + ... $ form a basis in
$L/zL$ and hence in $L$ considered as a free $\CC [z]$-module. Since 
$z\p J/\p \tau^{\d}\ \in zL \subset \L$, the 2nd derivatives 
$z\p^2 J /\p \tau^{\d}\p\tau^{\l}$ are in $L$ again and are uniquely
representable via the basis as $\sum A_{\d\l}^{\mu} \p J/\p \tau^{\mu}$
with the coefficients $A_{\d\l}^{\mu}(\tau) \in \CC [z]$. 
On the other hand, these 2nd derivatives are in $z\H_{-}$, and therefore 
the coefficients do not depend on $z$. 
Thus $\p J/\p \tau^{d}$ form a fundamental solution of the
pencil of flat connections depending linearly on $z^{-1}$:
\begin{equation} \label{Dmod} 
z \frac{\p}{\p \tau^{\l}} \frac{\p J}{\p\tau^{\d}} = 
\sum_{\mu} A_{\d\l}^{\mu} (\tau) \frac{\p J}{\p \tau^{\mu}} .\end{equation}
Obviously $A_{\d\l}^{\mu} = A_{\l\d}^{\mu}$. 
Moreover, the string flow $\f \mapsto \exp (u/z) \f $ preserves 
$-z+ z\H_{-}$ and projects along $\H_{-}$ to $-z+\tau \mapsto -z+\tau - u\1$.
This implies that $z\p J / \p \tau^{1} = J$ and respectively 
$(A_{1\l}^{\mu})$ is the identity matrix.
Furthermore, it is not hard to derive now that the multiplications
on the tangent spaces $T_{\tau}H = H$ given by 
\[ \phi_{\d}\bullet \phi_{\l} = \sum_{\mu} A_{\d\l}^{\mu}(\tau)\phi_{\mu} \]
define associative commutative algebra structures with the unit $\1$ 
which are Frobenius with respect to the inner product on $H$ 
(i.e. $(a\bullet b,c)=(a,b\bullet c)$) and satisfy the integrability condition
imposed by (\ref{Dmod}).  This is what a Frobenius structure
is (see for instance \cite{M}).\footnote{
To match this with Dubrovin's original definition \cite{D} one can use the 
generating function $\F$ to show that in fact 
$A_{\d\l}^{\mu}=\sum g^{\mu\nu} \ 
\lan \phi_{\nu} ,\phi_{\d},\phi_{\l}\ran (\tau)$ and therefore 
$\F (\tau,0,0,...)$ satisfies the WDVV-equation.} 
  
Conversely, given a Frobenius structure one recovers a $J$-function 
by looking for a fundamental solution to the system ({\em cf.} (\ref{Dmod}))
\[  z \frac{\p}{\p \tau^{\gl}} S = \phi_{\gl} \bullet S  \]
in the form of an operator-valued series
$S =\1+S_1(\tau)z^{-1}+S_2(\tau)z^{-2}+...$ satisfying $S^*(-z)S(z)=\1$.
Such a solution always exists and yields the corresponding $J$-function 
$ J^{\d}(z,\tau)=z [S^*(z,\tau)]^{\d}_1 $ 
and the Lagrangian cone $\L$ enveloping the family of Lagrangian spaces 
$L=S^{-1}(z,\tau) \H_{+}$ and satisfies the condition ($\star$). 
A choice of the fundamental
solution $S$ is called in \cite{GiQ} a {\em calibration} of the Frobenius
structure. The calibration is unique up to the right multiplication 
$S\mapsto S M$ by elements $M=\1+M_1z^{-1}+M_2z^{-2}+...$
of the ``lower-triangular'' subgroup in the twisted loop group. Thus the 
action of this subgroup on our class of cones $\L$ changes calibrations
but does not change Frobenius structures (while more general elements 
of $\L^{(2)}GL(H)$, generally speaking, change Frobenius structures
as well).

\medskip

As an example, consider the translation 
invariant Frobenius manifold 
defined by a given Frobenius algebra structure $\circ$ on $H$. 
The system (\ref{Dmod})
has constant coefficients and allows for the following obvious solution
\[ J(z,\tau) = z\ e^{(\tau \circ) /z } = 
\sum_{k\geq 0} z^{1-k}\frac{\tau^{\circ k}}{k!}  .\]
The corresponding cone is $\L=\{\ \e^{(\tau\circ)/z} \q(z),\ 
\q\in H[z]\ \}$. 
Infinitesimal automorphisms of $\L$ in the twisted loop 
Lie algebra contains the subalgebra
\[ z^{-1}(a_0\circ) + z (a_1\circ) + z^3(a_2\circ) +...,\ a_i\in H\]
and actually coincides with it when the algebra $(H,\circ)$ is 
semisimple.

\medskip

{\bf Quantum cohomology.} In Gromov -- Witten theory, the genus $0$ 
descendent potential $\F$ of a compact (almost) K\"ahler manifold $X$
is defined in terms on intersection theory on moduli spaces $X_{0,n,d}$
of degree $d\in H_2(X,\ZZ)$ stable maps $\gS \to X$ of genus $0$ 
complex curves with $n$ marked points \cite{K}. Namely  
\begin{equation} \label{untwisted} 
\F = \sum_{d,n} \frac{Q^d}{n!} \int_{[X_{0,n,d}]} 
\prod_{i=1}^n (\sum_{k=0}^{\infty} \ev_i^*(t_k)\psi_i^k ) .\end{equation}
Here $[X_{0,n,d}]$ is the virtual fundamental cycle \cite{LT}, 
$t_k\in H^*(X,\QQ)$ are cohomology classes of the target space $X$ to be 
pulled back to the moduli spaces by the evaluation maps  
$\ev_i: X_{0,n,d}\to X$ at the $i$-th marked points, 
and $\psi_i$ is the $1$-st Chern class of the {\em universal cotangent
line bundle} $L_i$. The fiber of $L_i$ is defined as the cotangent line 
$T^*_{\gs_i}\gS$ to the curve at the $i$-th marked point.

Thus the cohomology of the target space equipped with the
Poincar\'e intersection pairing 
\[ g_{\mu\nu} = \int_X \phi_{\mu}\phi_{\nu} \]
plays the role of the inner product super-space $H$. Yet the situation here
does not quite fit our axiomatic set-up since an appropriate
{\em Novikov ring} $\QQ [[Q]]$ takes the place of the ground field.
Moreover, one often has to consider the Lagrangian section 
$\L=\operatorname{graph} d\F$ in the completion of the space $\H$ consisting
of the series $\sum f_k z^k$ infinite in both directions but 
such that $f_k \to 0$ as $k\to +\infty$ in the $Q$-adic topology of 
$H=H^*(X,\QQ [[Q]])$. In what follows we will allow this and similar
generalizations of the framework without special notice.            
       
The genus $0$ descendent potential is known to satisfy the axioms 
DE, SE and TRR and therefore defines on $H$ a Frobenius structure
with the multiplication $\bullet$ known as {\em quantum cup-product}.
This Frobenius structure comes therefore equipped with the canonical
calibration $S$ and the $J$-function:  
\begin{equation} \label{J} J(z,\tau) = z+\tau+\sum_{d,n} \frac{Q^d}{n!} \ 
(\ev_{n+1})_* [\frac{\prod_{i=1}^n \ev_i^*(\tau)}{z-\psi_{n+1}}], 
\end{equation}
where $(\ev_{n+1})_*$ is the push-forward along the evaluation map
$\ev_{n+1}: X_{0,n+1,d} \to X$. The present description of this expression
as the defining section of the whole cone $\L$ seems more satisfactory than
the {\em ad hoc} construction in the literature on mirror formulas.

\medskip

Below we briefly recall several formulations from \cite{CG}
exploiting the language of Lagrangian cones in the context of genus $0$
Gromov -- Witten theory.

\medskip

{\bf Gravitational ancestors.} Forgetting the maps $\Sigma \to X$ 
(and maybe several last marked points) and stabilizing the remaining
marked curves one defines the {\em contraction maps} 
$\ct: X_{g,m+n,d}\to \M_{g,m}$
of moduli spaces of stable maps to the Deligne -- Mumford spaces.     
The latter carry their own universal cotangent line bundles and their
$1$-st Chern classes whose pull-backs to $X_{g,m+n,d}$ will be denoted 
$\bar{L}_i$ and $\bar{\psi}_i$ respectively. We call {\em ancestors}
the Gromov -- Witten invariants which ``mistakenly'' use the classes
$\bar{\psi}$ instead of $\psi$. A theorem of Kontsevich -- Manin \cite{KM}
expresses descendents via ancestors and vice versa. It was E. Getzler \cite{Ge}
who emphasized a primary role of ancestors in the theory of topological
recursion relations. The following theorem interprets their results
geometrically in the case $g=0$. 

Introduce the genus $0$ {\em ancestor potential} $\bar{\F}_{\tau}$ depending
on the parameter $\tau\in H$:
\[ \bar{\F}_{\tau} =  \sum_{d,m,n} \frac{Q^d}{m!n!} \int_{[X_{0,n,d}]} 
\prod_{i=1}^m (\sum_{k=0}^{\infty} \ev_i^*(\bar{t}_k)\bar{\psi}_i^k ) \ 
\prod_{j=m+1}^{m+n} \ev_j^*(\tau) .\]
Let $\L_{\tau}$ denotes the Lagrangian section in $\H=T^*\H_{+}$ 
representing the graph of $d\bar{\F}_{\tau}$ subject, as usual, to the
dilaton shift $\q(z) = -z + \sum \bar{t}_k z^k$.

\medskip

{\bf Theorem 2} (see Appendix $2$ in \cite{CG}). {\em We have
$ \L_{\tau} = S (\tau) \L$, where $\L$ represents the descendent potential,
and $S (\tau) \in \L^{(2)}GL(H)$ is the corresponding calibration.}

\medskip

The ancestor potential $\bar{\F}_0$ has zero $2$-jet when $\bar{t}_0=0$
since $\dim \M_{0,2+m} = m-1<m$. Therefore $\L_{\tau}$ is tangent to $\H_{+}$
along $z\H_{+}$. Thus Theorem $2$ implies directly (i.e. by-passing 
Theorem $1$) that the graph $\L$ of $d\F$ satisfies ($\star$).

\medskip

{\bf Corollary.} {\em The sections $\L_{\tau}$ are themselves Lagrangian
cones satisfying the condition ($\star$) 
(and isomorphic to each other and to $\L$).}

\medskip

The Frobenius structures defined by the cones $\L_{\tau}$ are isomorphic to
that for $\L$ since the calibration $S$ is ``lower-triangular''.

On the other hand, Theorem $2$ tells us how to define the ancestor 
potentials $\bar{\F}_{\tau}$ in the axiomatic theory: restrict the quadratic
form ${\mathbf p} (\q)/2$ to $S (\tau) \L$.
In particular, the genus $0$ ancestor potentials $\bar{\F}_{\tau}$
satisfy SE, DE and TRR (which agrees with some results in \cite{RT}).

Moreover, let us call a function $\G (\bar{t}_0,\bar{t}_1,...)$ {\em tame} if
\begin{equation} \label{tame} 
\frac{\p}{\p \bar{t}^{\a_1}_{k_1}}...\frac{\p}{\p \bar{t}^{\a_r}_{k_r}} 
\ \G\ |_{\bar{\t}=0} = 0\ \ \text{whenever}\ k_1+...+k_r > r-3 .\end{equation}
The ancestor potentials in Gromov -- Witten theory are tame on the dimensional
grounds: $\dim \M_{0,r} = r-3$. The following result shows that in the 
axiomatic theory there is no need to add this condition as a new axiom.

\medskip

{\bf Proposition.} {\em The ancestor potentials $\bar{\F}_{\tau}$
(defined in the axiomatic theory as the generating functions for the
Lagrangian cones $\L_{\tau}:=S (\tau) \L$) are tame.}

\medskip

{\em Proof.} The ancestor potentials themselves satisfy DE$+$SE$+$TRR
but in addition have zero $2$-jet when $\bar{t}_1=\bar{t}_{2}=...=0$.
Let us apply the derivation 
$\p^r /\p \bar{t}^{\a_1}_{n_1}...\p \bar{t}^{\a_r}_{n_r}$
to the identity TRR for $\bar{F}_{\tau}$ and assume that
(\ref{tame}) is satisfied for the terms on the right hand side. If 
at least one of these terms is non-zero at $\bar{\t}=0$, then we must have
$k+l+m+n_1+...+n_r \leq r-1$. Thus (\ref{tame}) is satisfied for the left
hand side. This allows us to derive (\ref{tame}) by induction on 
$k_1+...+k_r$. $\square$  
  
Apparently, a more general result applicable to higher genus theory is 
contained in \cite{Ge}.
 
\medskip

{\bf Twisted Gromov -- Witten invariants.}
Let $E$ be a holomorphic vector bundle over a K\"ahler target space $X$.
One can pull it back to a moduli space $X_{g,n+1,d}$ by the evaluation map
$\ev_{n+1}: X_{g,n+1,d}\to X$ and then push the result forward to 
$X_{g,n,d}$ by the map $\ft_{n+1}: X_{g,n+1,d}\to X_{g,n,d}$ forgetting
the $n+1$-st marked point. This way one gets an element in $K^0(X_{g,n,d})$:
\[ E_{g,n,d} = 
\ \text{a virtual bundle with the fibers}\ \ 
H^0(\Sigma; f^*E) \ominus H^1(\Sigma; f^*E).\] 
Fixing an invertible multiplicative characteristic class 
$S (\cdot ) = e^{\sum s_{k} ch_k (\cdot)}$
of complex vector bundles, one defines {\em twisted} GW-invariants
by systematically using intersection indices against the {\em twisted} 
virtual fundamental cycle $S(E_{g,n,d})\cap [X_{g,n,d}]$. 
When $E$ is a negative line bundle, and $S$ is
the inverse total Chern class, the twisted theory represents GW-invariants
of the total space of the bundle. When $E$ is a positive line bundle,
and $S$ is the Euler class, the twisted genus $0$ theory is very close
to the GW-theory of the hypersurface defined by a section of the line
bundle. The following ``Quantum Riemann -- Roch theorem'' 
describes the twisted version $\F^{tw}$ of the genus $0$ descendent potential 
(\ref{untwisted}) via the untwisted one.

Given $E$ and $S$, represent the the twisted genus $0$ descendent potential
by the Lagrangian submanifold $\L^{tw} \subset \H$ obtained from
the graph of $d\F^{tw}$ using the dilaton shift convention 
$\q (z) = \sqrt{S(E)} (\t(z)-z)$. 

\medskip

{\bf Theorem 3} (see \cite{CG}). 
{\em The Lagrangian submanifold $\L^{tw}$ 
is obtained from the Lagrangian cone $\L$ of the untwisted theory by a
linear symplectic transformation:
\[ \L^{tw} = \triangle \L, \ \ 
\triangle \sim \sqrt{S(E)} \prod_{m=1}^{\infty} S(E\otimes L^{-m}) \]
where $L$ is a line bundle with $c_1(L)=z$.}

\medskip

More precisely, put $s(x)=\sum s_k x^k/k!$ and let $x_i$ be the Chern roots of 
$E$, so that $S(E)=\exp \sum_i s(x_i)$. Then the expression 
\[ \ln \{ S(E)^{1/2} \prod_{m=1}^{\infty} S(E\otimes L^{-m}) \} = 
\sum_i [ \frac{s(x_i)}{2} +\sum_{m=1}^{\infty} s(x_i-mz) ] \]
has the asymptotical expansion 
\[  \frac{1}{2}\sum_i \frac{1+e^{-z\p/\p x_i}}{1-e^{-z\p/\p x_i}}\ s(x_i) =
  \sum_{0\leq m}\ \ \sum_{0\leq l\leq \dim_{\CC}X}
s_{2m-1+l} \frac{B_{2m}}{(2m)!} ch_l(E) z^{2m-1},\] 
where $B_{2m}$ are the Bernoulli numbers, and $ch_l(E)$ are understood
as operators of multiplication in the cohomology algebra $H$ of the target 
space $X$. The operator $\ln \triangle$ is the infinitesimal
symplectic transformation defined by this series. 

\medskip

When expressed in the original terms of the generating functions,
the relationship stated in the theorem would have the form
\[ \F^{tw}(\q) = \ \text{critical value in $\q'$ of}\ \ 
\G(\q,\q') + \F (\q'),\]
where $\G$ is the quadratic generating function of the above linear 
symplectic transformation. The operation of computing the critical value
is complicated,  and it seems not obvious if the
new function $\F^{tw}$ satisfies, say, the TRR. In fact it is not hard to
see directly that the twisted genus $0$ GW-theory satisfies the same axioms
DE, SE and TRR. Then another question arises: are these equations for $\F^{tw}$
formal consequences of those for $\F$ or they carry some additional information
about the untwisted theory? The use of Lagrangian cones makes the answer 
obvious: since the linear transformation in the theorem belongs to the 
group $\L^{(2)}GL(H)$, the condition ($\star$) for
$\L^{tw}$ is formally equivalent to that for $\L$. This elementary 
example illustrates possible advantages of the geometrical viewpoint.

\medskip 

{\bf Quantum Lefschetz Theorem.} The following result is 
a formal (but combinatorially non-trivial) consequence of Theorem $3$.
 
Note that a Lagrangian cone satisfying ($\star$) is determined by its generic 
$\dim H$-dimensional submanifold. We begin with the $J$-function $J(-z,\tau)$ 
representing the Lagrangian cone $\L$ of the untwisted theory on $X$ and 
generate a new Lagrangian cone $\L'$ satisfying ($\star$) 
by using the following
function $I(-z,\tau)$ in the role of the ``generic $\dim H$-dimensional
submanifold'':     
\[ I(z,\tau) \sim \int_0^{\infty} e^{x/z} J(z, \tau+\rho \ln x)\ dx \ 
\text{\Huge /}\ \int_0^{\infty} e^{(x-\rho \ln x)/z} dx ,\]
where $\rho =\gl+c_1(E)$ is the total Chern class 
of a complex line bundle $E$ over $X$, 
and the integrals should be understood in the sense of stationary phase
asymptotics as $z\to 0$.    

\medskip

{\bf Theorem 4} (see \cite{CG}). {\em The Lagrangian cone $\L'$
determined by the function $I$ coincides with the cone $\L^{tw}$ of
the genus $0$ Gromov -- Witten theory on $X$ twisted by the line bundle $E$
and by the total Chern class
\[  S(\cdot)=\gl^{\dim (\cdot)}+
c_1(\cdot )\gl^{\dim (\cdot )-1}+...+c_{\dim(\cdot)}(\cdot). \]}

\medskip

This is an abstract formulation of a very general result which in the
limit $\gl\to 0$ allows one to compute genus $0$ GW-invariants of
a hypersurface in $X$ defined by a section of $E$ in a form
generalizing Quantum Lefschetz Theorems of 
B. Kim, A. Bertram, Y.-P. Lee and A. Gathmann 
\cite{BCKS, Ber, Lee1, Gat}.
Namely, Theorem $4$ applies beyond ``small'' quantum cohomology theory,
to general type complete intersections as well, and also offers a new insight 
on the nature of {\em mirror maps} as a very special case of {\em Birkhoff 
factorization} in loop groups. 
 In particular, given the $J$-function of
a toric Fano manifold $X$, the theorem yields the mirror formulas of
\cite{Gi2} for the $J$-function of a Fano or Calabi -- Yau 
complete intersection in $X$ (including the celebrated example of
quintic three-folds in $\CC P^4$).\footnote{
Ironically, in order to find the $J$-functions of general Fano toric manifolds,
one still has to apply the methods of fixed point localization used in
the original proof \cite{Gi3, Gi2} of the mirror formulas for toric
complete intersections.}

\medskip

{\bf Singularity theory.} Although not named so, Frobenius structures first 
emerged around 1980 from K. Saito's theory of primitive forms in 
singularity theory (see \cite{H} for a modern exposition).

Let $f:(\CC^m,0)\to (\CC, 0)$ be the germ of a holomorphic function 
(for simplicity --- weighted homogeneous) at an isolated critical point of 
multiplicity $N$. Let $F(x,\tau)$ be a miniversal deformation of $f$, i.e. 
a family of functions in $\x\in \CC^m$
depending on the $N$-dimensional parameter $\tau\in \Tau$ and such that 
$F(x,0)=f(x)$
and $\p F/\p \tau^{\a}|_{\tau=0},\ \a=1,...,N$, represent a basis in the 
local algebra $H$ of the critical point:
\[ H  = \CC [x] / (f_{x_1},...,f_{x_m}).\]
The tangent spaces $T_{\tau}\Tau$ are canonically identified with the algebras
$\CC [x] / (F_x)$ of functions on the critical sets. To make the
algebras $T_{\tau}\Tau$ Frobenius, one picks a holomorphic volume form
$\go $ on $\CC^m$ (possibly depending on $\tau$) and introduces the 
{\em residue pairing} on $T_{\tau}\Tau$ via
\[  (a , b)  = (\frac{1}{2\pi i})^m \oint ... \oint
 a(y)\ b(y)\ \frac{dy_1\w ...\w dy_m}{F_{y_1}\ ...\ F_{y_m}} .\]
(We assume here that $y$ is a unimodular coordinate system on $\CC^m$,
i.e. $\go = dy_1\w ... \w dy_m$.)
To make $\Tau$ a Frobenius manifold one takes the volume form 
$\go $ to be {\em primitive}, i.e. satisfying very special conditions
which in particular guarantee that the residue metric is flat. 

Many attributes of abstract Frobenius structures make more sense when
interpreted in terms of singularity theory. For example, the system 
(\ref{Dmod})  
is satisfied by the complex oscillating integrals
\[  \J_{\Gb}(\tau) = (- 2\pi  z)^{-m/2}\int_{\Gb} e^{F(x,\tau)/z} \go \]
written in flat coordinates $\tau$ of the residue metric on $\Tau$.

The oscillating integrals also satisfy the following weighted
homogeneity condition:
\begin{equation} \label{E}
(z\frac{\p}{\p z} + \sum (\deg{\tau^{\l}}) \tau^{\l} 
\frac{\p}{\p \tau^{\l}})
\ z\frac{\p}{\p \tau^{\d}} \J = -\mu_{\d}\ z\frac{\p}{\p \tau^{\d}} \J ,
\end{equation}
where $\mu_{\d} =  \deg (\p_{\d} F) +\deg (\go)- m/2,\ \d=1,...,N$, is the 
{\em spectrum} of the singularity symmetric about $0$. The 
pencil (\ref{Dmod}) of flat connections over $\Tau$ can therefore be
extended in the $z$-direction by the operators
\[ \nabla_{\tau}:= \frac{\p}{\p z} +\frac{\mu}{z} +(E\bullet)/z^2. \]
Here $E=\sum (\deg \tau^{\l}) \tau^{\l} \p/\p \tau^{\l}$ is the 
{\em Euler field} and $\mu =\operatorname{diag}(\mu_1,...,\mu_N)$ 
is the {\em Hodge grading operator} (anti-symmetric with respect 
to the metric and diagonal in a graded basis of $H$). 

The extended connection is flat (since the oscillating integrals
$z\p \J_{\Gb} /\p \tau^{\d}$ provide a basis of flat sections)
and can be considered as an {\em isomonodromic} family of connections
$\nabla_{\tau}$ over $z\in \CC\backslash 0$ depending on the parameter 
$\tau\in \Tau$. In this situation, the calibration 
$S (z,\tau) =\1+S_1z^{-1}+S_2z^{-2}+...$
can be chosen as a gauge transformation near $z=\infty$
of the operator $\nabla_{\tau}$ to the normal form $\nabla_0=\p /\p z +\mu/z$. 
In other words, the fundamental solution matrix 
$(z\p \J_{\Gb_{\l}}/\p \tau^{\d})$
consisting of oscillating integrals can be written as 
$S (z,\tau)\ z^{-\mu}\ C$
(where $C$ is a constant matrix), and the $J$-function can
be extracted from $\J$ in a similar way.   

In singularity theory, gravitational descendents are not defined 
(at least mathematically) in an intrinsic way, but rather recovered
from the axioms, i.e. the total descendent potential is defined as 
the generating function for the Lagrangian cone $\L \subset \H$ 
determined by the $J$-function.   

\medskip

{\bf Semisimplicity.} In the case $\dim H=1$ (interpreted either
as the Gromov -- Witten theory of the point target space or in mirror
terms as $A_1$ singularity) the $J$-function $z e^{\tau/z}$
generates the descendent potential (in dilaton-shifted variables)
\begin{equation} \label{A_1}
\F (\q) = \ \text{the critical value of}\ \frac{1}{2}\int_0^{\tau} 
(q_0+q_1x+...+q_k\frac{x^k}{k!}+...)^2 \ dx. \end{equation}
Consider now the model Cartesian product case of a semisimple 
translation invariant Frobenius structure defined by the {\em algebra}
$\CC^N$. Let $\L^{(N)}$ be the corresponding Lagrangian submanifold.
It is isomorphic to $(\L^{(1)})^N$ and is generated by direct sum 
$\F (\q^{(1)})+...+\F (\q^{(N)})$ of $N$ copies of the function (\ref{A_1}).

\medskip

{\bf Theorem 5.} {\em The Lagrangian submanifold $\L \subset \H$
representing the germ of a semisimple $N$-dimensional Frobenius manifold
is locally isomorphic to $\L^{(N)} \subset \CC ^N ((z^{-1}))$.}

\medskip

{\em Proof.} This is a reformulation of a result from \cite{GiS} 
({\em cf.} Exercise $3.7$ in \cite{D} though) about 
existence of an asymptotical fundamental solution  
\[ (z \frac{\p J^{\l}}{\p \tau^{\d}}) 
\sim \Psi(\tau) R (z,\tau) e^{U(\tau)/z} \]
to the system (\ref{Dmod}) at semisimple points $\tau$. 
Here $R =\1+R_1z+R_2z^2+... $ is an ``upper-triangular'' element of
the (completed) twisted loop group $\L^{(2)}GL_N$, 
$\Psi: \CC^N \to H$ is an isomorphism
of the inner product spaces, and $U=\operatorname{diag} (u_1,...,u_N)$ is
the diagonal matrix of Dubrovin's {\em canonical coordinates} \cite{D} on
the Frobenius manifold. 

In the case of K. Saito's Frobenius structures of 
singularity theory, $u_{\l}(\tau)$ are critical values of the Morse functions
$F(\cdot, \tau)$, and an asymptotical solution $\Psi R \exp (U/z)$
can be found as the stationary phase asymptotics of the complex 
oscillating integrals $z\p \J_{\Gb_{\l}}/\p \tau^{\d}$. 
   
In general, consider the cones $\L_{\tau}=S (\tau) \L$ representing the 
ancestor potentials. The operators 
$\exp (-U(\tau)/z)R^{-1}(\tau)\Psi^{-1}(\tau)$ act on the cones $\L_{\tau}$ 
despite the fact that $R=\1+R_1z+...$ are infinite $z$-series. This is
because the ancestor potentials are tame and in particular can be 
considered as series in $\bar{t}_0,\bar{t}_1$ with coefficients 
{\em polynomial} in $\bar{t}_2,\bar{t}_3,...$ (see Section $8$ in \cite{GiA}
for more details). The resulting cone $\L'$ does not 
depend on $\tau$ since the operators $S$ and $\Psi R \exp (U/z)$ 
satisfy the same differential equations (\ref{Dmod}) in $\tau$. 
On the other hand, the operators $R$ preserve the isotropic space 
$z\CC^N [[z]]$ contained in (the completion of) the cones $\L_{\tau}$.    
This implies that the cone $\L'$ together with the point $-z$ contains
all points $\exp (-U/z)(-z)$. Thus $z \exp (U/z)$ is the $J$-function of $\L'$
which shows that $\L'=\L^{(N)}$. $\square$ 

\medskip

The operator $R (z,\tau)$ is unique up to the right multiplication by 
the automorphisms $\exp (a_0/z+a_1z+a_2z^3+...),\ a_i\in \CC^N$, of 
$\L^{(N)}$\footnote
{and up to reversing or renumbering coordinates in $\CC^N$.}.
The ambiguity in is eliminated by imposing the additional homogeneity
condition (\ref{E}) (or equivalently $L_E R =0$) available in the presence 
of the Euler vector field $E$ in the definition of the Frobenius structure.
This is the case in singularity theory and in (non-equivariant!) 
Gromov -- Witten theory.    

\medskip

{\bf Virasoro constraints.} The so called ``Virasoro conjecture''
was invented by T. Eguchi, K. Hori, M. Jinzenji, C.-S. Xiong and S. Katz
\cite{EHX, EJX} and upgraded to the rank of an axiom in abstract 
topological field theory by B. Dubrovin and Y. Zhang \cite{DZ1}. 
In the axiomatic context, the Virasoro constraints are defined
in the presence of the additional grading axiom expressing the role
of the Euler field $E$ in the definition of Frobenius structures. 
The genus $0$ Virasoro constraints are known to contain no information
in addition to the grading axiom and the system DE$+$SE$+$TRR. 
The geometrical argument below replaces
the original proof of this fact given by X. Liu -- G. Tian 
\cite{LT} as well as the shorter proof given by E. Getzler \cite{Ge1}.

Suppose that linear operators $A$ and $B$ on a vector space are anti-symmetric
with respect to a bilinear form. Then the operators 
$l_m = ABABA ... BA$ ($B$ repeated $m$ times) are also anti-symmetric.
On the other hand, if $AB-BA = B$, then $l_{m}, m=0,1,2, ...$, commute as
the vector fields $x^{m+1}\p/\p x$ on the line: $[l_m,l_n]=(n-m)l_{m+n}$.

In the symplectic space $(\H, \gO)$, consider infinitesimal symplectic
transformations $A,B$ satisfying $[A,B]=B$, where 
$A=l_0=z d/dz + 1/2 + a(z)$ with $a^*(-z)+a(z)=0$, and  $B$ is 
multiplication by $z$. 

The additional grading axiom for the descendent potential $\F$ is formulated as
the invariance of the cone $\L$ under the flow of the 
linear hamiltonian vector field on $\H$ defined by the operator $l_0$ with 
a special choice of $a$ (determined by the the Euler field $E$).  For instance,
in the case of weighted - homogeneous singularities $a=\mu$, and in 
Gromov -- Witten theory take  $a = \mu+\rho /z$ where $\mu$ is the Hodge 
grading operator on $H^*(X;\QQ)$ and $\rho$ is the operator of multiplication 
by $c_1(T_X)$ in $H^*(X;\QQ)$. However the explicit form of $a$ is irrelevant
for the present discussion. Note that the vector field defined
by the operator $l_{-1}=B^{-1}=z^{-1}$ (corresponding to $\p /\p x$) is 
also tangent to $\L$ due to the string equation. 
The following theorem expresses the genus $0$ Virasoro constraints for $\F$.

\medskip

{\bf Theorem 6.} {\em Suppose the vector field 
on $\H$ defined by the operator $l_0$ is tangent to the Lagrangian cone
$\L$ satisfying the condition ($\star$). Then the same is true 
for the vector fields defined by the operators 
$l_m=l_0zl_0z...zl_0,\ m=1,2,...$.}

\medskip

{\em Proof.} 
Let $L$ be the Lagrangian space tangent to $\L$ along $zL$. We know therefore
that for any $\f \in zL$ the vector $l_0\f$ is in $L$. This implies 
$zl_0\f \in zL$, implies $l_0zl_0\f \in L$, implies $zl_0zl_0\f \in zL$,
..., and therefore for all $m=1,2,...$ the vectors $l_m\f$ are tangent to 
the cone at the point $\f$. $\square $.

\medskip  

{\bf Higher genus and quantization.} Correlators of the higher genus axiomatic
theory are arranged into the Taylor coefficients of a sequence of
formal functions $\F^{(g)}(\t)$, $g=0,1,2,...$,
called {\em genus $g$ descendent potentials}, which includes the genus $0$
theory with $\F = \F^{(0)}$. In this section, we collect some results from
\cite{GiS, GiQ, CG, GiA, GiM, JaK, JoK} illustrating the following general 
observation:
{\em the higher genus theory is quantization of the genus $0$ theory 
and inherits the twisted loop group of symmetries}.

The polarization $\H=\H_{+}\oplus \H_{-} =T^*\H_{+}$ allows one to quantize
quadratic hamiltonians on the symplectic space $(\H,\gO)$ using the following
standard rules:
\[ (q_{\a}q_{\b})\hat{\ } = q_{\a}q_{\b} /\h, \ 
   (q_{\a}p_{\b})\hat{\ } = q_{\a}\p/\p q_{\b},\ 
   (p_{\a}p_{\b})\hat{\ } = \h \p^2/\p q_{\a} \p q_{\b} .\]
Quantized quadratic hamiltonians act on functions of $\q \in \H_{+}$
depending on the parameter $\h$.
They form a {\em projective} representation of the Lie algebra of
infinitesimal symplectic transformations and of its Lie subalgebra 
$\L^{(2)}gl (H)$. The notation $\hat{A}$ for a simplectic transformation
$A\in \L^{(2)}GL (H)$ will be understood as the operator $\exp (\ln A)
\hat{\ }$.

The functions $\F^{(g)}$ are assembled into the {\em total descendent 
potential} 
\[ \D = e^{ \F^{(0)}/\h + \F^{(1)} + \h \F^{(2)} + \h^2 \F^{(3)} + ...} \]
Expressions of this form where $\F^{(g)}$ are formal functions of $\t$
will be called {\em asymptotical functions} and, after the dilaton shift
$\q (z) = -z + t_0+t_1z+...$,  considered as ``asymptotical elements of
the Fock space''. This means that we treat the
$1$-dimensional subspace $\lan \D \ran$ spanned by $\D$ as a point 
in the projective space --- the potential domain of our quantization
representation. 

\medskip

\noindent $\bullet$ Theorem $3$ above is derived in \cite{CG} as the 
quasi-classical limit $\h \to 0$ of the higher genus 
``Quantum Riemann -- Roch Theorem'' expressing twisted Gromov -- Witten
invariants via the untwisted ones:  
\[ \lan\D^{tw}\ran = \hat{\triangle} \lan \D \ran \]
where $\triangle$ is the symplectic transformation in Theorem $3$.  
This result is based on Mumford's Grothendieck -- Riemann -- Roch 
formula \cite{Mu} applied to the universal families of stable maps
and generalizes the result of Faber -- Pandharipande \cite{FP} for Hodge
integrals corresponding to the case of the trivial twisting bundle $E$.

\medskip  
 
\noindent $\bullet$ The higher genus Kontsevich -- Manin
formula \cite{KM} (see also \cite{Ge}) 
relating descendents and ancestors in Gromov -- Witten theory has been 
rewritten in \cite{GiQ} (see also \cite{CG}) as the quantization of the
calibration operators $S(\tau)$ of Theorem $2$:   
\[ \A_{\tau} = \hat{S} (\tau) \D ,\]
where $\A_{\tau} = \exp \sum \h^{g-1} \bar{\F}_{\tau}^{(g)}(\bar{\t})$ 
is the total ancestor potential.

\medskip  
 
\noindent $\bullet$ According to \cite{GiS, GiQ}, the machinery of fixed point
localization and summation over graphs used by
M. Kontsevich \cite{K} for computation of equivariant Gromov -- Witten 
invariants in the case of tori acting on the target spaces with 
isolated fixed points gives rise to the following formula for the
equivariant total ancestor potentials:
\begin{equation} \label{R}
 \lan \A_{\tau} \ran = \Psi (\tau) \ \hat{R} (\tau)\ \exp (U/z)\hat{\ }
\ \lan \D_{pt}^{\otimes N} \ran .\end{equation}
Here the rightmost function is the product $\D_{pt} (\q^{(1)}) ... 
\D_{pt}^{(N)} $ of $N=\dim H$ copies of the total descendent potential
$\D_{pt}$ of the one-point space $X=pt$, and $\Psi R \exp (U/z)$ is the
asymptotical solution to the system (\ref{Dmod}).\footnote
{See \cite{GiS,GiQ} for the explanation of how to
fix the ambiguity in the construction of the asymptotical solution
in the equivariant quantum cohomology theory lacking the Euler vector
field $E$.}

\medskip  
 
\noindent $\bullet$ The formula (\ref{R}) makes sense for any {\em semisimple}
Frobenius structure and thus can be used as a definition of the higher
genus ancestor potentials in the abstract theory. According to \cite{GiQ} 
this definition meets a number of expectations. Namely: 

$\circ $ The total ancestor potential $\A_{\tau}$ thus defined is {\em tame}
in the sense that the corresponding functions $\bar{\F}^{(g)}_{\tau}(\bar{\t})$
satisfy
\[ \frac{\p}{\p \bar{t}^{\a_1}_{k_1}}...\frac{\p}{\p \bar{t}^{\a_r}_{k_r}} 
\ \bar{\F}^{(g)}_{\tau}\ |_{\bar{\t}=0} = 0\ \ \text{whenever}\ 
k_1+...+k_r > r+3g-3 .\]

$\circ$ The corresponding total descendent potential {\em defined} by
$\lan \D \ran =  \hat{S}^{-1} (\tau)  \lan \A_{\tau} \ran $ does not 
depend on $\tau$.\footnote{Together with the previous property this means that
$\D$ satisfies the $3g-2$-jet condition  of Eguchi -- Xiong \cite{EX,Ge1} 
playing the role of TRR in the higher genus theory.}

$\circ$ The string equation $\hat{l}_{-1} \D = 0$ always holds true,
and the higher Virasoro constraints $\hat{l}_m \D =0, m=1,2,...$ hold true 
if the grading condition $\hat{l}_0 \D \in \lan \D \ran$ is satisfied.

\medskip  
 
\noindent $\bullet$ In non-equivariant Gromov--Witten theory, the formula 
\begin{equation} \label{D} \lan \D \ran = \hat{S}^{-1}\ \Psi\ \hat{R}\ 
e^{(U/z)\hat{\ }} \ \lan \D_{pt}^{\otimes N} \ran \end{equation} 
is known to be correct for complex projective spaces and other 
Fano toric manifolds \cite{GiQ}
and for complete flag manifolds \cite{JoK} (which in particular
implies that the Virasoro conjecture holds true for these manifolds).
Also, it is an easy exercise on application of the dilaton equation for 
$\D_{pt}$ to check that the result of Jarvis -- Kimura \cite{JaK} 
computing gravitational descendents in the 
Gromov -- Witten theory on the quotient {\em orbifold} of the one - point 
space by a finite group agrees with (\ref{D}). In this example, the Frobenius
structure is translation invariant, $(H, \circ, (\cdot,\cdot ))$ 
is the center of the group ring, $\Psi$ describes its orthonormal
diagonalization while $S = R =\1$.

\medskip  
 
\noindent $\bullet$ According to the part of Witten's conjectures \cite{W}
proved by M. Kontsevich \cite{KA} (see also \cite{OP})
the asymptotical function $\D_{pt}$
describing intersection theory on the Deligne -- Mumford spaces $\M_{g,m}$
satisfies the KdV hierarchy of evolution equations with the time variables
$t_0,t_1,t_2,...$. The remaining part of the conjectures similarly identifies
the total descendent potential in the moduli theory of curves
equipped with {\em $n$-spin structures} with solutions to the 
$n$KdV (or {\em Gelfand -- Dickey}) hierarchies of integrable systems.
Corresponding Frobenius structure coincides \cite{W} with K. Saito's 
structure on the miniversal deformation of the $A_{n-1}$-singularity. 
It is a result of \cite{GiA} that in the case of $A_{n-1}$-singularity
the function (\ref{D}) indeed satisfies the equations of the $n$KdV hierarchy.
In \cite{GiM} this result is generalized to the $ADE$-singularities.
 
\medskip

Conventional wisdom says that the structure of the axiomatic 
higher genus theory should reflect the geometry of the Deligne -- Mumford 
spaces, and vice versa. Yet both subjects are far from being clear. 
In view of the above examples it is tempting to impose the following 
two requirements on conjectural axioms of the higher genus theory. 

{\em A. The set of asymptotical elements of the Fock space satisfying
the axioms has to be invariant under the quantization representation
of the twisted loop Lie algebra.}

In particular, the stabilizer of a Lagrangian cone $\L$ with the generating
function $\F$ satisfying the axioms of the genus $0$ theory has to act on 
the set of asymptotical elements satisfying the axioms of the higher genus 
theory and having the same genus $0$ part $\F^{(0)}=\F$. 

{\em B. The action of the stabilizer is transitive.}

The latter requirement is inspired by the uniqueness lemma of 
Dubrovin -- Zhang \cite{DZ} which says that in the {\em semisimple} case
an asymptotical function $\D$ satisfying the higher genus TRR and the 
Virasoro constraints is unique (and thus coincides with (\ref{D}) according 
to the results quoted earlier). It would be preferable however 
to derive the Virasoro constraints from other axioms (such as or 
similar to A$+$B) just the way it is done in the proof of Theorem $6$ in
the genus $0$ theory.
 
\medskip

{\bf Quantum K-theory.} The material of this section represents our  
joint work with Tom Coates and shows that the description of the correlators 
and their properties in geometrical terms of Lagrangian cones and quantization 
remains valid in the context of quantum K-theory. 

\medskip

In quantum K-theory, a basis of ``observables'' has the form
$\Phi_{\a} L_i^k$ where $\{ \Phi^{\a} \}$ is a linear $\QQ$-basis 
in $K^*(X)$ (here $K^0(X)$ is the Grothendieck group of complex vector 
bundles over the target space $X$) and $L_1,...,L_n$ 
are the universal cotangent line 
bundles over the moduli spaces of stable maps $X_{g,n,d}$. 

Correlators
can be defined in cohomological terms using the {\em virtual tangent bundles} 
$T_{g,n,d}$ of the moduli spaces $X_{g,n,d}$ and with holomorphic Euler 
characteristics in mind:
\begin{equation} \label{chi} 
\lan \Phi^{\a_1} L^{k_1}, ... , \Phi^{\a_n} L^{k_n}\ran_{g,n,d} = 
\int_{[X_{g,n,d}]} \td (T_{g,n,d}) \ \ch_* [ 
\otimes_{i=1}^n \ev_i^*(\Phi^{\a_i}) \otimes L_i^{\otimes k_i} ]  .
\end{equation}

First properties of K-theoretic Gromov -- Witten invariants are discussed in
\cite{GiK}, and the foundations in the setting of algebraic target spaces have
been laid down by Yuan-Pin Lee in \cite{Lee2}. 
We have to stress that the correlators 
(\ref{chi}) are only approximations to the actual holomorphic Euler 
characteristics which in the orbifold / orbibundle context are given
by Kawasaki's Hirzebruch -- Riemann -- Roch Theorem \cite{Ka}. The
correlators (\ref{chi}) differ therefore from those used in the papers  
\cite{GiK,Lee2}. However the general {\em properties} of the correlators
(as opposed to their {\em values}) remain the same as in \cite{GiK, Lee2}.
This is because the bundles, sheaves and their properties used in those 
papers are the same as in the ``fake'' version of K-theory
considered here. 

In complete analogy with the cohomology theory, one can introduce K-theoretic 
descendent potentials $\F^{(g)}(\t)$ as formal functions 
of $t_0+t_1L+t_2L^2+...$ with $t_i \in K:=K^*(X)\otimes \QQ [[Q]]$.
As it was found in \cite{GiK}, K-theoretic genus $0$ Gromov -- Witten 
invariants
define a ``Frobenius-like'' structure on $K$ which however exhibits 
the following remarkable distinction from the
case of cohomology theory. The constant coefficient metric
\[ g_{\a\b} = \lan \Phi_{\a}, \1 , \Phi_{\b} \ran_{0,3,0} = 
\chi (X; \Phi_{\a} \otimes \Phi_{\b}) = \int_{X} \td (T_X)\ \ch_*(\Phi_{\a})\
\ch_*(\Phi_{\b}) \]
has to be replaced in all formulas by the non-constant one:
\[ G_{\a\b}(\tau) = g_{\a\b}+\sum_{n\geq 0,d} \frac{Q^d}{n!} 
\lan \Phi_{\a}, \tau, ..., \tau, \Phi_{\b} \ran_{0,n,d},\ \ 
\text{where}\ \tau =\sum \tau^{\mu}\Phi_{\mu} \in K .\]   
For instance, the K-theory version of the $J$-function (\ref{J})  
\begin{equation} \label{JK} 
J = (1-q) + \tau + \sum_{n\geq 0,d} \frac{Q^d}{n!} (\ev_{n+1})^K_* 
\left[ \frac{\otimes_{i=1}^n \ev_i^*(\tau)}{\1 - q L_{n+1}} \right] 
\end{equation}
satisfies the system 
\begin{equation} \label{DK}
(1-q) \frac{\p}{\p \tau^{\d}} \frac{\p}{\p \tau^{\l}} J = 
\sum  A_{\d\l}^{\mu} \ \frac{\p}{\p \tau^{\mu}} \ J. \end{equation}
The multiplication 
$\Phi_{\d}\bullet \Phi_{\l}=\sum_{\mu} A_{\d\l}^{\mu}\Phi_{\mu}$ 
makes the tangent spaces $T_{\tau}K$ Frobenius algebras with respect to
the inner product $(G_{\mu\nu})$ since in fact  
\[ A_{\d\l}^{\mu}(\tau) = \sum_{\nu} G^{\mu\nu}(\tau) \ 
\p_{\d,0}\p_{\l,0}\p_{\nu,0} \F^{(0)}|_{\t =(\tau,0,...)}, \ \ 
\text{where}\ (G^{\mu\nu})=(G_{\mu\nu})^{-1} .\]
Moreover, the system (\ref{DK}) at $q=-1$ yields, as it is shown in \cite{GiK},
the Levi -- Civita connection of the metric $(G_{\mu\nu})$ which is therefore
flat. 

\medskip

Is it possible to adapt our language of symplectic loop spaces to
absorb these, rather dramatic, changes? 
We explaine below that the answer is positive 
at least in the version of quantum K-theory discussed here.\footnote
{In the original version of quantum K-theory studied in \cite{GiK,Lee2}
the orbifold features of the moduli spaces $X_{g,n,d}$ play 
prominent role. For example, the classes $L_i \in K^0(X_{g,n,d})$ no longer
behave unipotently. This causes some difficulties which seem
formal at first (so that an analogue of Theorem $7$ below still holds
true), but become overwhelming in the more sophisticated  
Riemann -- Roch problems discussed in the next section.}  
  
\medskip

First, note that the $J$-function (\ref{JK}) considered as a $1/(1-q)$-series
becomes finite when reduced modulo $Q^d$ since $L_{n+1}-1$ is nilpotent
in $K^0(X_{0,n,d})$. We will call such Laurent $1/(1-q)$-series 
{\em convergent} away from $q=1$ in the $Q$-adic topology 
(i.e. for those $q\in \QQ [[Q]]$ whose $Q$-adic distance to $1$ is $\geq 1$). 

Next, consider the loop space $\K$ of all vector Laurent series with 
coefficients in $K$ convergent in this sense away from $q=1$.   
Equip $K$ with the inner product $(\Phi_{\a}, \Phi_{\b})=g_{\a\b}$
(the original constant one!) and define in $\K$ the symplectic form
\[ \gO (\f , \g) = \frac{1}{2\pi i} \oint (\f(q^{-1}), \g (q))\ \frac{dq}{q} 
:= - [ \Res_{q=0}+\Res_{q=\infty} ]  (\f (q^{-1}),\g (q))\ \frac{dq}{q} .\]   
Substituting $q=\exp z$ (as motivated by $\ch_*(L_{n+1})=\exp \psi_{n+1}$) 
we find that this $\gO$ is in a sense the same as in cohomology theory. 
Consider the polarization $\K = \K_{+}\oplus \K_{-}$
where $\K_{+} = K [q]$ and $\K_{-}$ consists of the series convergent away
from $q=1$ {\em and vanishing at $q=\infty$}. The polarization is Lagrangian 
(although not invariant with respect to $z\mapsto -z$) and 
identifies $\K$ with (a topologized version of) $T^*\K_{+}$. 
We use the dilaton shift $\q (q) = 1-q+t_0+t_1q+t_2q^2+... $ to
encode the genus $0$ descendent potential by a Lagrangian
submanifolds $\L$, and furthermore --- 
to identify the total descendent and ancestor potentials of higher
genus quantum K-theory with the appropriate asymptotical
functions $\D$ near $\q = (1-q) \in \K_{+}$.
 
Now, repeat the construction for ancestors $\L_{\tau}$, $\A_{\tau}$
using the inner product $(\Phi_{\a},\Phi_{\b})_{\tau} = G_{\a\b}(\tau)$
and the corresponding symplectic structure $\gO^{\tau}$ in the space
$\K^{\tau}=\K$. 

It turns out that in quantum K-theory the calibration operators $S (q,\tau)$
defined via the $J$-function (\ref{JK}) as
\[ S^{\a}_{\b}(q,\tau) = \sum_{\mu\nu} G^{\a\mu}(\tau)\ 
\frac{\p J^{\nu}}{\p \tau^{\mu}} (q,\tau)\ g_{\nu\b} \] 
satisfy 
\[ \sum_{\a\b} S^{\a}_{\mu} (q,\tau) G^{\a\b}(\tau) S_{\nu}^{\b}(q^{-1},\tau) 
= g_{\mu\nu}. \]
This shows that $S (\tau)$ is symplectic as a linear map from 
$(\K,\gO)$ to $(\K^{\tau},\gO^{\tau})$. 
The following result is the K-theoretic version of the 
Kontsevich -- Manin theorem \cite{KM} relating descendents and ancestors.

\medskip

{\bf Theorem 7}. {\em We have $e^{F^{(1)}(\tau)} \A_{\tau} = \hat{S} (\tau) \D$
and in particular $\L_{\tau} = S (\tau) \L$.}

\medskip

Since 
$\ch_*(\bar{L}_i-1)=e^{\bar{\psi}_i}-1$ are nilpotent in $H^*(\M_{0,n};\QQ)$, the 
Lagrangian sections $\L_{\tau}$ are tangent to $\K_{+}^{\tau}$ along
$(1-q)\K_{+}^{\tau}$.

\medskip

{\bf Corollary.} {\em The Lagrangian submanifold $\L\subset \K$ is a cone 
with the vertex at the origin and satisfies the condition that its
tangent spaces $L$ are tangent to $\L$ exactly along $(1-q)L$.}

\medskip

The conic property of $\L$ is equivalent to the genus $0$ case of the
K-theoretic {\em dilaton equation} 
\[ \lan \t^{(1)}(L), ..., \t^{(n)}(L), 1-L\ran_{g,n+1,d} = (2-2g-n) 
\lan \t^{(1)}(L), ..., \t^{(n)}(L)\ran_{g,n,d} \ran ,\]
where  $\t^{(i)}$ are arbitrary polynomials of the universal cotangent line 
bundles $L=L_i$. 

The operator of multiplication by  
\[ \frac{1}{1-q}-\frac{1}{2}=\frac{1}{2}\ \frac{1+q}{1-q} \] 
is anti-symmetric with respect to $\gO$ and defines a linear hamiltonian 
vector field on $(\K,\gO)$ which is tangent to the cone $\L$. This
property of $\L$ expresses the following genus $0$ K-theoretic 
{\em string equation} 
\begin{align} \notag 
\lan \t^{(1)}(L), ..., \t^{(n)}(L), 1\ran_{0,n+1,d} &\ = & \\ 
\notag \lan \t^{(1)}(L), ..., \t^{(n)}(L)\ran_{0,n,d}\ + & \sum_{i=1}^n &
\lan \t^{(1)}(L), ..., \frac{\t^{(i)}(L)-\t^{(i)}(1)}{L-1},..., 
\t^{(n)}(L)\ran_{0,n,d}. \end{align}

\medskip

As we mentioned earlier, the same results hold true in the original version
of the quantum K-theory studied in \cite{GiK, Lee2}.

\medskip

{\bf Quantum cobordism theory.} The complex cobordism theory $MU^*(\cdot)$
is defined in terms of homotopy classes of maps to the spectrum 
$MU(k)$ of the Thom spaces of universal $U_{k/2}$-bundles:
\[ MU^n(B) = \lim_{k\mapsto \infty} \pi (\Sigma^k B, MU(n+k)) .\]
When $B$ is a stably almost complex manifold of real dimension $m$, 
the famous Pontryagin -- Thom construction identifies elements of 
$MU^n(B)$ with 
appropriately framed bordism classes of maps $M \to B$ of stably almost 
complex manifolds $M$ of real dimension $m-n$. This identification
plays the role of the Poincar\'e isomorphism.  
Similarly to the complex K-theory, 
there is the {\em Chern -- Dold character} which provides natural 
multiplicative isomorphisms
\[ \Ch: MU^*(M)\otimes \QQ \to H^*(M, \Lambda^*). \]
Here $\Lambda^*=MU^*(pt)\otimes \QQ$ is the coefficient ring of the
theory and is isomorphic to the polynomial algebra on the generators 
of degrees $-2k$ 
Poincar\'e -- dual to the bordism calsses $[\CC P^k]$. 
The theory $MU^*$ is known to be the universal extraordinary cohomology 
theory where complex vector bundles are oriented.
Orientation of complex bundles in $MU^*$
is uniquely determined by the cobordism-valued Euler class 
$u\in MU^2(\CC P^{\infty})$ of the
universal complex line bundle. 
Explicitly, the Euler class of ${\mathcal O} (1)$ over $\CC P^N$ is 
Poincar\'e -- dual to the embedding  
$\CC P^{N-1}\to \CC P^N $ of a hyperplane section.
The image of $u$ under the Chern -- Dold character has the form
\begin{equation} \label{u} u (z) = z + a_1 z^2 + a_2 z^3 + ... \end{equation}
where $z$ is the cohomological $1$-st Chern of the universal line bundle
${\mathcal O} (1)$, and $\{ a_k \} $ is another set of generators in 
$\Lambda^*$.
The operation of tensor product of line bundles with the Euler classes
$v $ and $w$ defines a formal group law $F(v,w)$ on 
$MU^*(\CC P^{\infty}) = \Lambda^*[[u]]$. The 
series $u(z)$ is interpreted as an isomorphism with the additive group
$(x,y) \mapsto x+y$: $F(v,w) = u( z(v) + z(w) )$.   
Here $z(\cdot)$ is the series inverse to $u(z)$. It is known as 
the {\em logarithm} of the formal group and explicitly
takes on the form 
\[ z = u + [\CC P^1]\frac{u^2}{2}+[\CC P^2]\frac{u^3}{3} + 
[\CC P^3]\frac{u^4}{4} + ... \]
Specialization of the parameters $[\CC P^k] \mapsto 0$ yields the 
cohomology theory, and $[\CC P^k] \mapsto 1$ yields the complex K-theory.
In the latter example, $z=-\ln (1-u)$ and hence $u(z) = 1-\exp (-z)$.
Similarly to the K-theory, one can compute push-forwards in $MU^*$-theory
in terms of cohomology theory. In particular, for a stably almost complex
manifold $B$, we have the {\em Hirzebruch -- Riemann -- Roch formula}
\begin{equation} \label{HRR} \forall c\in MU^*(B),\ \  \pi_*^{MU}(c) = 
\int_{B} \Ch (c)\ \Td (T_B) \ \in  
\ \Lambda^*,\end{equation}  
where $\pi : B \to pt$, and $\Td (T_B)$ is the {\em Todd genus} 
of the tangent bundle. 
It is characterized as the only multiplicative 
characteristic class which for the universal line bundle is equal to
\[ 
 \Td = \frac{z}{u(z)} = \exp \sum_{k=1}^{\infty} s_k\frac{z^k}{k!} .
\] 
Here $s_1,s_2,...$ is one more set of generators in $\Lambda^*$. 
To round up the introduction, let us mention the Lanweber -- Novikov 
algebra of stable cohomological operations in complex cobordism theory.
The operations correspond to cobordism-valued characteristic classes $\sigma$
of complex vector bundles. 
To apply such an operation to the cobordism element of $B$ 
Poincar\'e-dual to a map $\pi: M\to B$ of stably almost complex manifolds,
one takes the relative normal bundle $\pi^*T_B\ominus T_M$ over $M$ 
and pushes-forward its characteristic class $\sigma$ to $B$.
According to Buchshtaber -- Shokurov \cite{BS}, after tensoring
with $\QQ$ the Landweber -- Novikov algebra becomes isomorphic to the 
algebra of left-invariant differential
operators on the group of diffeomorphisms (\ref{u}).
The Landweber -- Novikov operations commute with the Chern -- Dold 
character and can therefore be expressed in the cohomology theory 
$H^*(\cdot,\Lambda^*)$ as certain differential operators
on the algebra $\Lambda^*$ of functions on the group. The algebra is generated
by the derivations $L_n$ whose action on the generators
$a_k$ is given by
\[ (L_n u)(z) = u(z)^{n+1},\ \ \text{or, equivalently,}\ \ 
L_n = u^{n+1}\p /\p u, \ n=1,2,... .\]   
These generators correspond to the characteristic classes defined by the
Newton polynomials 
\[ u_1^n+u_2^n+... \in MU^*(\CC P^{\infty}\times \CC P^{\infty}\times ... ). \]

\medskip

In a sense, the idea of Gromov -- Witten invariants with values
in cobordisms is already present in Gromov's original philosophy \cite{Gr} 
of symplectic invariants to be constructed as bordism invariants of spaces of 
pseudo-holomorphic curves.  The possibility to define 
Gromov -- Witten invariants with values in the cobordism ring $\Lambda^*$ 
is mentioned in Kontsevich's work \cite{K}. This proposal 
was further advanced by Morava \cite{Mo} in a hope to explain the Virasoro
constraints in terms of the Landweber -- Novikov operations. Following
suggestions of these authors, we define cobordism-valued Gromov -- Witten
invariants via (\ref{HRR}) using the virtual tangent bundles 
$T_{g,n,d}$ of the moduli spaces $X_{g,n,d}$. We formulate below 
the genus $0$ version of the ``Quantum Hirzebruch -- Riemann -- Roch 
Theorem'' 
which expresses these new Gromov -- Witten invariants via the old ones.
The discussion below represents joint work of Tom Coates and the author.
The forthcoming thesis \cite{C} contains many details omitted here. 

\medskip

In quantum cobordism theory, a basis of observables has the form
$\Phi_{\a} u_i^k$ where $\{ \Phi^{\a} \}$ is a basis of 
$H\Lambda :=MU^*(X)\otimes \QQ [[Q]]$ over $\Lambda^*\otimes \QQ[[Q]]$, 
and $u_i$ is the Euler class of the universal cotangent line bundle $L_i$.
The correlators are defined by the formula
\begin{equation} \label{corr}  
\lan \Phi_{\a_1} u^{k_1}, ... , \Phi_{\a_n} u^{k_n}\ran_{g,n,d} := 
\int_{[X_{g,n,d}]} e^{\sum_{k=1}^{\infty} s_k \ch_k (T_{g,n,d})} \ 
\prod_{i=1}^n \left[ \ev_i^* \Ch (\Phi_{\a_i})\ u(\psi_i)^{k_i} \right] .
\end{equation} 
One can introduce $\Lambda^*$-valued descendent potentials $\F^{(g)}(\t)$
in complete analogy with the cohomology theory as formal functions 
of $t_0+t_1u+t_2u^2+...$ with $t_i \in H\Lambda$.

In order to express the cobordism-valued Gromov -- Witten invariants 
in terms of cohomological ones one needs to control the
classes $\ch_k(T_{g,n,d})$. Roughly speaking, the virtual tangent
bundles $T_{g,n,d}$ consist of two parts. One of them, representing 
deformations of holomorphic {\em maps} of Riemann surfaces with a
fixed complex structure, coincides with $E_{g,n,d}$ (as in Theorem $3$ about
twisted Gromov -- Witten invariants) with $E = T_X$. The effect of this
part can be described therefore via Theorem $3$. The other part 
represents deformations of complex structures and contributes into the 
correaltors in a complicated fashion. Surprisingly, the very formalism of 
the symplectic loop space takes effective care of these contributions.     

\medskip

We denote by $\Lambda = \QQ \{ s_1,s_2,... \}$ the coefficient ring of 
complex cobordism theory completed with respect to the $\s$-adic 
norm defined by the grading $\deg s_k = -k$.
Let $H\Lambda $ be the cobordism group of the target space $X$ over $\Lambda$
equipped with the inner product 
\[ (a,b)_{\s} :=\pi_*^{MU}(ab) = \int_X \Ch (a) \Ch (b) 
\exp \left( \sum s_k \ch_k(T_X) \right).\]
Denote $\H\Lambda$ the loop space $H\Lambda \{ \{ u^{-1} \} \}$ of 
Laurent series $\sum f_k u^k$ with coefficients $f_k\in \H\Lambda$
which can be non-zero for all $k\in \ZZ$ but should satisfy the condition
that $\f_k \to 0$ in the $\s$-adic topology as $k\to +\infty$. 
Introduce the symplectic form with values in $\Lambda $: 
\[ \gO_{\s} (\f,\g) :=  
\frac{1}{2\pi i} \oint \pi_*^{MU}(\f (u^*)\ \g(u))\ \sum_{k\geq 0} 
[\CC P^k] u^k\ du = \frac{1}{2\pi i} \oint (\f(u(-z),\g(u(z))_{\s}\ dz .\]
The following {\em quantum Chern -- Dold character} identifies the symplectic
structure $\gO_{\s}$ with its cohomological version $\gO =\gO_{\0}$:
\[ \qCh: \H\Lambda \to \H =\lambda\otimes H ((z^{-1})): \ \ \qCh (\f) = 
\sqrt{\Td (T_X)} \ \sum_{k\in \ZZ} \Ch (f_k) u^k(z).\]
Assuming $|\psi |<|z|$ write
\[ \frac{1}{u(-z-\psi)} = \sum_{k\geq 0} u(\psi)^k v_k(u(z)) = 
\sum_{k\geq 0} \psi^k w_k(z).\]
Put 
\[ \H\Lambda_{+}=H\Lambda \{ u \} =\ \{ \ \text{power series $\sum q_ku^k$ 
with $q_k\to 0$} \ \} \] and 
\[ \H\Lambda_{-}=\ \{ \ \text{arbitrary infinite series $\sum p_k v_k(u)$}\ \}
.\]
The following residue computation shows that $\{ ..., u^k, ..., v_k(u),...\}
\ $\ (as well as\ 
$\{ ..., z^k, ..., w_k(z),... \}$) 
is a Darboux basis in $(\H\Lambda, \gO_{\s})$ and implies that the spaces
$\H\Lambda_{\pm}$ form a Lagrangian polarization (depending on $\s$):
\[ \frac{1}{2\pi i} \oint \frac{dz}{u(z-x)\ u(-z-y)} = 
\left\{ \begin{array}{ccl} 0 & \text{if} & |z|<|x|,|y|\ \text{or}
\ |x|,|y|<|z|\\
\pm \frac{1}{u(-x-y)}& \text{if} &
|x|<|z|<|y| \ \text{or}\ |y|<|z|<|x| \end{array} \right. . \]  
Introduce the dilaton shift convention
\[ \q (u) = \t (u) + u^*(u),\ \ \text{where}\ \ u^*(u(z)) = u(-z).\]
We consider the genus $0$ descendent potential $\F^{(0)}$ of the quantum
cobordism theory as a function depending formally on the parameters 
$\s=(s_1,s_2,...)$ and define the corresponding family of Lagrangian 
submanifolds $\L_{\s}$ in the family of symplectic spaces 
$T^*\H\Lambda_{+} = \H\Lambda $ over $Spec \L$. The following
theorem expresses $\L_{\s}$ in terms of the cone $\L = \L_{\0}$ describing
the quantum cohomology theory of $X$.
      
\medskip 

{\bf Theorem 8} (see \cite{C}). {\em The image of $\L_{\s}$ under the
quantum Chern -- Dold character coincides with the Lagrangian cone of 
the Gromov -- Witten theory on $X$ twisted by the class $\Td (T_X)$:
\[ \qCh (\L_{\s}) = \square \L, \ \ 
\square \sim \sqrt{\Td (T_X)} \ \prod_{m=1}^{\infty} \Td (T_X \otimes L^{-m}),
\]
where $L$ is a line bundle with $c_1(L)=z$.}

\medskip

In fact $\qCh (\H\Lambda_{+}) = \H_{+}$, but $\qCh (\H\Lambda_{-})$ and
the dilaton shift depend on $\s$. This causes a discrepancy between the 
descendent potentials of the quantum cobordism theory and the twisted theory. 
Remarkably, the discrepancy accounts for the entire contribution of 
variations of complex structures on Riemann surfaces
into the virtual tangent bundles of the moduli spaces of stable maps.

The same happens in the higher genus theory. The quantum Chern -- Dold 
character identifies the Heisenberg Lie algebras of the spaces 
$(\H\Lambda, \gO_{\s})$. By the virtue of the Stone -- von Neumann theorem
and the Schur lemma, the Fock spaces corresponding to different polarizations
are identified projectively. The total descendent potentials $\D_{\s}$ of 
the quantum cobordism theory differ from the total descendent potentials
$\D_{\s}^{tw}$ of the appropriately twisted Gromov -- Witten theory 
in such a way that the corresponding lines $\lan \D_{\s}\ran$ and 
$\lan \D_{\s}^{tw}\ran $ in the representation space of the Heisenberg 
algebra coincide.

\medskip

Returning to the genus $0$ case, notice that in the case of translation
invariant Frobenius structure the cone $\L$ is invariant under 
the symplectic transformations defined by $\square $.

\medskip

{\bf Corollary.} {\em When $X=pt$, we have $\qCh (\L_{\s}) = \L$.}

\medskip

The formal group law
\begin{equation}  \label{lambda} u(x+y)=u(x)+u(y)-\lambda u(x)u(y),\ \ 
u(z) = (1-e^{-\lambda z})/\lambda, \end{equation}
interpolates between cohomology and K-theory. Thus the specialization
$[\CC P^k] \mapsto \lambda^k$ allows one to adjust the quantum 
Hisrzebruch -- Riemann -- Roch theorem  
to the case of quantum K-theory.  (The correlators discussed in the 
previous section correspond to $\lambda = -1$.)

\medskip

{\bf Corollary.} {\em The quantum Chern character 
$\qch: \K \mapsto \H$ defined via $\qch u = (1-\exp (-\lambda z))/\lambda$
identifies $\qch(\L_{\lambda})$ with $\square \L$ where 
$\square$ is defined via the Todd class 
$\td (L) = \lambda z/(1-\exp (-\lambda z)$.} 

\medskip

The quantum Chern -- Dold character transforms multiplication by 
$u$ to multiplication by $u(z)=z+... \in \Lambda \{ z\}$. 
Transformations defined by $\square$ belong to the twisted 
loop group and thus commute with multiplication by $u(z)$. 

\medskip

{\bf Corollary.} {\em The submanifolds $\L_{\s} \subset \H\Lambda$
are Lagrangian cones satisfying the condition ($\star$).}

\medskip

In particular, the dilaton equation holds true in quantum cobordism theory
(which is easy to prove directly for any genus). Moreover, the cubical
form on $L/zL$ defined by the correlators 
$\lan \Phi_{\a},\Phi_{\b},\Phi_{\c}\ran $ represents the {\em Yukawa coupling}
defined on any Lagrangian submanifold in a linear symplectic space
(see \cite{GiH}). The form coincides therefore with the structure tensor 
of the Frobenius manifold defined by the twisted Gromov -- Witten theory.
However in quantum cobordism theory, there is no a simple 
formula for the string equation, as there is no general reasons for
flatness of the metric $G_{\a\b}=\lan \Phi_{\a}, 1, \Phi_{\b}\ran$ or 
associativity of the quantum cup-product whose constructions 
depend on the polarization. In fact, a key step ---
Barannikov's derivation of the equation (\ref{Dmod}) for the $J$-function ---
was based on the relation $u^*\H\Lambda_{-}=H\Lambda+\H\Lambda_{-}$.
It is not hard to see that the group laws (\ref{lambda}) are the only ones
satisfying this condition. 

Indeed, the inclusion $\forall k\geq 0, \ H u^*w_k\subset 
H\Lambda+\H\Lambda_{-}$ means that the projection of all $u(-z)w_k(z)$ 
along the subspace $Span (w_l(z), l=0,1,2,...)$ yields constant polynomials.
By definition we have $\sum_l w_l(-z)x^l = 1/u(z-x)$ when 
$|x|<|z|$ and $\sum_k u(-z)w_k(z)y^k = u(-z)/u(-z-y)$ when $|y|<|z|$.
In terms of these generating functions the inclusion is equivalent 
therefore to the condition that
\[ \frac{1}{2\pi i} \oint_{|x|,|y|<|z|} \frac{u(-z) \ dz}{u(z-x)\ u(-z-y)}
= \frac{u(-x)}{u(-x-y)} -\frac{u(y)}{u(-x-y)} \]
is independent of $x$ for all $y$. Differenting 
$u(-x-y) f(y) = u(-x)-u(y)$ in $x$ we find $u'(-x-y) f(y) = u'(-x)$ which
at $x=0$ yields $f(y)=1/u'(-y)$ (since $u'(0)=1$) and implies 
$u'(-x-y)=u'(-x) u'(-y)$. Thus $u'(-x)=\exp (\gl x)$ and respectively 
$u(z)=(1-\exp (-\gl z))/\gl$.  
 
It would be interesting to find out if ellipticity of $u (z)$ brings 
any ``improvements'' in properties of {\em elliptic quantum cohomology}
in comparison with the general quantum cobordism theory. 


\enddocument

\bibitem{Ar} V. I. Arnold. 

\bibitem{AVG} V. I. Arnold, S. M.  Gusein-Zade, A. N. Varchenko. 
{\em Singularities of differentiable maps. Vol. II. Monodromy and asymptotics of integrals.} 
Monographs in Mathematics, 83. Birkh\"auser Boston, Inc., Boston, MA, 1988. viii+492 pp

\bibitem{CG} T. Coates, A. Givental. {\em Quantum Riemann -- Roch, Lefschetz and Serre.}
\newline arXiv: math.AG/0110142.

\bibitem{D} B. Dubrovin. {\em Geometry of 2D topological filed theories.} In:
Integrable Systems and Quantum Groups. 
Springer Lecture Notes in Math. 1620 (1996), 120--348. 

\bibitem{DZ} B. Dubrovin, Y. Zhang. {\em Normal forms of hierarchies of integrable PDEs, 
Frobenius manifolds and Gromov -- Witten invariants.} arXive: math.DG/0108160.


\bibitem{GiQ} A. Givental. {\em Gromov -- Witten invariants
and quantization of quadratic hamiltonians.} Moscow Mathematical Journal, v.1(2001), no. 4, 
551--568.

\bibitem{GiS} A. Givental. {\em Semisimple Frobenius structures at higher genus.} 
Intern. Math. Res. Notices, 2001, No. 23, 1265--1286. 

\bibitem{GiE} A. Givental. {\em Elliptic Gromov -- Witten invariants and
the generalized mirror conjecture.} In: Integrable Systems and Algebraic
Geometry. World Sci. Publ., River Edge, NJ, 1998, 107--155.
 
\bibitem{H} C. Hertling. {\em Frobenius manifolds and moduli spaces for singularities.}
Cambridge Tracts in Mathematics. Cambridge University Press, 2002, 280 pp. 

\bibitem{K} V. Kac. {\em Infinite dimensional Lie algebras.} 3rd edition. 
Cambridge University Press, 1990, 400 pp.

\bibitem{Ko}  M. Kontsevich {\em Intersection theory on the moduli space of
curves and the matrix Airy function.} Commun. Math. Phys. {\bf 147} 
(1992), 1 -- 23.

\bibitem{Ma} Yu. I. Manin. {\em Frobenius manifolds, quantum cohomology,
and moduli spaces.} AMS Colloquium Publ. 47, Providence, RI, 1999, 303 pp.

\bibitem{PV} A. Polishchuk, A. Vaintrob. {\em Algebraic construction of Witten's top Chern class.}
arXiv: math.AG/0011032.

\bibitem{S} K. Saito. {On a linear structure of the quotient variety 
by a finite reflection group.} Publ. Res. Inst. Math. Sci. 29 (1993),
no. 4, 535--579. 

\bibitem{Sch} A. Schwarz. {\em On some mathematical problems of 2D-gravity and
$W_{h}$-gravity.} Modern Physics Letters A, Vol. 6, No. 7 (1991), 611 -- 616.

\bibitem{SW} G. Segal, G. Wilson. {\em Loop groups and equations of KdV type.}
 Publ. IHES, No. 61 (1985), 5--65.

\bibitem{Va} A.N. Varchenko, {\em Local residue and the intersection form in vanishing cohomology.} 
Izv. Akad. Nauk SSSR Ser. Mat. 49 (1985), no. 1, 32--54; English translation in:
Math. USSR Izvestiya, v. 26, no. 1 (1986), 31 -- 52.

\bibitem{W} E. Witten. {\em Two-dimensional gravity and intersection theory
on moduli space.} Surveys in Diff. Geom.  1 (1991), 243--310.

\end{thebibliography}

\enddocument